\documentclass[11pt]{amsart}

\usepackage[usenames,dvipsnames,svgnames,table]{xcolor}
\usepackage[hyphens]{url}
\usepackage[pagebackref,linktocpage=true,colorlinks=true,
   linkcolor=Blue,citecolor=BrickRed,urlcolor=RoyalBlue]{hyperref}
\usepackage[msc-links,abbrev]{amsrefs}
\usepackage{amsmath,amsthm,amssymb}
\usepackage{mathtools}
\usepackage{enumitem}

\newtheorem{theorem}{Theorem}[section]
\newtheorem{proposition}[theorem]{Proposition}
\newtheorem{lemma}[theorem]{Lemma}
\newtheorem{corollary}[theorem]{Corollary}
\theoremstyle{definition}
\newtheorem{note}[theorem]{Note}

\newcommand{\R}{\mathbf R}
\newcommand{\Z}{\mathbf Z}
\newcommand{\Sph}{\mathbf S}
\newcommand{\B}{\mathbf B}
\newcommand{\vol}{\operatorname{vol}}
\newcommand{\ip}[2]{\left\langle #1,#2\right\rangle}
\newcommand{\norm}[1]{\left| #1\right|}

\title[Shades and Convexity of Hypersurfaces]
 {Shades and Convexity of Hypersurfaces}

\author{Mohammad Ghomi}
\address{School of Mathematics, Georgia Institute of Technology,
   Atlanta, GA 30332}
\email{ghomi@math.gatech.edu}
\urladdr{ghomi.math.gatech.edu}

\date{\today\,(Last Typeset)}


\begin{document}

\begin{abstract}
Let $f\colon M^n\to\R^{n+1}$, $n\geq2$, be a smooth immersion of a closed
oriented manifold, with Gauss map $\nu\colon M\to\Sph^n$. For each unit
vector $u$, the \emph{shade} $S_u\subset M$ is
the set of points where $\nu$ makes an acute angle with $u$.
We show that $f$ is a convex embedding if and only if every shade is
contractible. More generally, if every shade is $k$-connected, then
$f$ is a convex embedding provided that $n\leq 2k+1$, or $n=2k+2$ and
$H_{k+1}(M;\Z)=0$. 
On the other hand, for every $n\geq 2k+3$ there is a nonconvex embedding
$f\colon\Sph^n\to\R^{n+1}$ all of whose shades are $k$-connected.
These results give sharp answers to a question of H. Wente in all dimensions.
\end{abstract}

\maketitle

\section{Introduction}

Let $M^n$, $n\geq2$, be a closed oriented manifold, and $f\colon M\to\R^{n+1}$ be a smooth immersion with Gauss map
$\nu\colon M\to\Sph^n$. For $u\in\Sph^n$, the \emph{shade} \cite{ghomi2007} (a.k.a.\ shadow \cite{ghomi2002-solution}) of $M$ in
direction $u$ is given by
\[
S_u\coloneqq\big\{p\in M: \ip{\nu(p)}{u}>0\big\}.
\]
If $f$ is a \emph{convex embedding}, i.e., $f(M)$ bounds a convex set, then $S_u$ is the region of $M$ which lies in the dark when $M$ is
illuminated by light rays parallel to $u$. In particular $S_u$ is connected. Motivated by stability problems for surfaces of constant mean curvature \cite{ghomi2002-solution, choe1990}, H. Wente asked in 1978 whether connectedness of the shades forces $f$ to be a convex embedding.
The author showed \cites{ghomi2002-shadows} that for $n=2$ the answer is yes
if the shades are simply connected, or else $M\cong\Sph^2$, i.e., $M$ is homeomorphic to $\Sph^2$; otherwise there are counterexamples. Here we extend those results to
higher dimensions. A space is \emph{$k$-connected} if  it is path connected and its homotopy groups $\pi_i=0$ for $1\leq i\leq k$.

\begin{theorem}\label{thm:main}
The immersion $f$ is a convex embedding if and only
if every shade is contractible.
More generally, if every shade is $k$-connected,
then $f$ is a convex embedding, provided that
$n\leq 2k+1$, or $n=2k+2$ and $H_{k+1}(M;\Z)=0$.
\end{theorem}

When $k\leq1$ and every shade is $k$-connected, the hypothesis
$H_{k+1}(M;\Z)=0$ for $n=2k+2$ is equivalent to $M\cong\Sph^n$;
see Note~\ref{note:borderline}.
So if $n\leq3$ or $M\cong\Sph^4$, it suffices that every shade be
simply connected; and if $M\cong\Sph^2$, it suffices that every shade
be connected. The dimension restrictions in Theorem~\ref{thm:main}
are sharp:

\begin{theorem}\label{thm:counterexamples}
For every $k\geq0$ and $n\geq 2k+3$, there is a smooth
nonconvex embedding $f\colon\Sph^n\to\R^{n+1}$
such that every shade is $k$-connected.
\end{theorem}

Thus, when $M\cong\Sph^n$, $k$-connected shades force $f$ to be a convex
embedding precisely when $n\leq 2k+2$. The proof of Theorem~\ref{thm:main}
is based on the equality case
of the Chern--Lashof theorem.
A lifting principle for projections of open sets
translates contractibility of the shades
into a count of the preimages of regular values of the Gauss map,
which yields convexity.
When the shades are only assumed to be $k$-connected,
we first show that they are contractible,
via Alexander duality.
The examples of
Theorem~\ref{thm:counterexamples}
are certain embeddings of the
join $\Sph^p*\Sph^q$ in $\R^{p+1}\oplus\R^{q+1}$
which are invariant under rotations of both factors.

\section{Convexity Characterizations}\label{sec:dim3}
Throughout we assume that $M$ is a \emph{closed} manifold, i.e., compact, connected without boundary, of dimension $n\geq 2$.
In this section we establish Theorem~\ref{thm:main}. The proof is
based on the Chern--Lashof theorem
\cites{chern-lashof1957,chern-lashof1958,kuiper1970,kuiper:tight},
which states that the \emph{total absolute curvature}
\[
\tau(f)\coloneqq\int_M|\det d\nu|
\]
satisfies $\tau(f)\geq\vol(\Sph^n)$,
with equality if and only if $f$ is a convex embedding.
By Sard's theorem and the area formula,
$\tau(f)=\int_{\Sph^n}\#\nu^{-1}(u)$, which is equal to $\vol(\Sph^n)$
when every regular value of $\nu$ has exactly one preimage.
So we obtain the following criterion.

\begin{lemma}\label{lem:chern-lashof}
If every regular value of $\nu$
has exactly one preimage,
then $f$ is a convex embedding.
\end{lemma}

\subsection{Contractible shades}
First we show that contractibility of the shades forces every
regular value of $\nu$ to have exactly one preimage.
To this end we write
$S_u=\{\sigma_u>0\}$, where $\sigma_u\coloneqq\ip{\nu}{u}$ is the
\emph{shading function}. Let
$$
E\coloneqq\big\{(p,u)\in M\times\Sph^n:\sigma_u(p)>0\big\},
$$
 and
$\pi\colon E\to\Sph^n$ be the projection into the second factor.
Then $\pi^{-1}(u)=S_u\times\{u\}$, so the shades of $M$ form the
fibers of $\pi$.
Furthermore, $\pi$ is surjective, i.e., every shade is nonempty.
Indeed $\sigma_u$ is not identically zero, since $\nu=\pm u$ at a
maximum point of the height function $\ip{f}{u}$, and so
$\sigma_u=\pm1\neq 0$ there. On the other hand
$\int_M\sigma_u=0$ by Stokes' theorem, since $\sigma_u$ times the
volume form of $M$ is the pullback under $f$ of the $n$-form
$\Omega(u,\cdot,\dots,\cdot)$ on $\R^{n+1}$, where $\Omega$ is the
volume form; this $n$-form is closed, and hence exact.
So $\sigma_u$ changes sign.
Although $\pi$ need not be a fiber bundle, the
openness of $E$ still admits enough local lifting to relate the
topology of $S_u$ to that of $E$:

\begin{lemma}\label{lem:micro}
Let $X$ and $B$ be manifolds, $E\subset X\times B$ be open, and
$\pi\colon E\to B$ be the second projection. Assume $\pi$ is surjective.
If the fibers of $\pi$ are $k$-connected, then $\pi$ induces bijections
on $\pi_i$ for $0\leq i\leq k$. If the fibers are contractible, then
$\pi$ is a homotopy equivalence.
\end{lemma}

\begin{proof}
Let $\B^m$ denote the closed unit ball in $\R^m$.
Given a lifting problem over $[0,1]\times \B^m$ with a lift at time zero,
openness of $E$ allows one to keep the $X$-coordinate fixed for a short
time. Compactness of $\B^m$ gives a uniform positive time. Thus $\pi$ is a
Serre microfibration. Recall that a map is \emph{$m$-connected} if it
induces bijections on $\pi_i$ for $i<m$, and a surjection on $\pi_m$.
If the fibers of $\pi$ are $k$-connected, then each of them is mapped
to a point by a $(k+1)$-connected map. So \cite[Thm.~1.3]{raptis2017},
applied with $q$ the identity map of $B$, shows that $\pi$ is
$(k+1)$-connected, which yields the first assertion. If the fibers are
contractible, then $\pi$ is a Serre fibration by a lemma of Weiss
\cite[Lem.~2.2]{weiss2005}, and so is a weak homotopy equivalence,
which is a homotopy equivalence by Whitehead's theorem,
since manifolds have the homotopy type of CW complexes \cite{milnor1959}.
\end{proof}

For $u\in\Sph^n$, let
$
H_u\coloneqq\{x\in\Sph^n:\ip{x}{u}>0\}
$
be the open hemisphere centered at $u$.
Then $S_u=\nu^{-1}(H_u)$,
and the fiber of $E$ over each $p\in M$ is $H_{\nu(p)}$.
We also put $I_u\coloneqq E\cap(M\times H_u)$.

\begin{proposition}\label{lem:count}
If every shade is contractible, then every regular value of
$\nu$ has exactly one preimage.
\end{proposition}

\begin{proof}
Let $-u$ be a regular value of $\nu$, so that $\nu^{-1}(-u)$ is finite,
and put $X\coloneqq M\setminus\nu^{-1}(-u)$.
Projection of $I_u$ to $H_u$ has contractible fibers $S_v$.
So Lemma~\ref{lem:micro} makes it a homotopy equivalence, and since
$H_u$ is contractible, so is $I_u$.
Projection of $I_u$ to $M$, on the other hand, has fiber
$H_u\cap H_{\nu(p)}$ over each $p\in M$, which is geodesically convex,
and is nonempty exactly when $\nu(p)\neq-u$.
So the image of this projection is $X$, and its fibers are
contractible. Thus Lemma~\ref{lem:micro} shows that $X$ is
contractible as well.
Write $N\coloneqq\#\nu^{-1}(-u)$. Since $X$ is contractible,
the long exact sequence of $(M,X)$ yields
$H_n(M,X;\Z)\cong H_n(M;\Z)\cong\Z$. On the other hand,
$H_n(M,X;\Z)\cong\Z^N$ by excision. Thus $N=1$.
\end{proof}

\begin{corollary}\label{cor:contractible}
$f$ is a convex embedding if and only if every shade is contractible.
\end{corollary}

\begin{proof}
If every shade is contractible, then $f$ is a convex embedding by
Proposition~\ref{lem:count} and Lemma~\ref{lem:chern-lashof}.
Conversely, each shade of a convex embedding is a graph over an open
convex set, and is therefore contractible; see \cite{ghomi2002-shadows}.
\end{proof}

\subsection{\texorpdfstring{$k$}{k}-connected shades}
Now we show that $k$-connected shades are contractible under the
hypotheses of Theorem~\ref{thm:main}; so Corollary~\ref{cor:contractible}
applies. We need three preliminary lemmas. For $u\in\Sph^n$ and $0<c<1$, we
define the \emph{enlarged shades}
$$
S_{u,c}\coloneqq\{\sigma_u>-c\}.
$$

\begin{lemma}\label{lem:enlarged}
If every shade is $k$-connected, then so is every enlarged shade.
\end{lemma}

\begin{proof}
Let $D\subset\Sph^n$ be the open cap of angular radius
$\arcsin c$ centered at $u$, and put
$Y\coloneqq E\cap(M\times D)$.
Projection of $Y$ to $D$ has fiber $S_v$ over each $v\in D$.
Since $D$ is contractible, Lemma~\ref{lem:micro} shows that
$Y$ is $k$-connected.
Projection of $Y$ to $M$, on the other hand, has fiber
$D\cap H_{\nu(p)}$ over each $p\in M$, which is geodesically convex,
and is nonempty exactly when the angle between $\nu(p)$ and $u$
is less than $\pi/2+\arcsin c$, i.e., when $\sigma_u(p)>-c$.
So the image of this projection is $S_{u,c}$, and its fibers
are contractible. Thus Lemma~\ref{lem:micro} shows that
$S_{u,c}$ is $k$-connected as well.
\end{proof}

\begin{lemma}\label{lem:homology}
If every shade is $k$-connected, then $M$ is
$\min\{k,n-1\}$-connected. Furthermore, if $n\leq 2k+2$, then
$H_i(M;\Z)=0$ for $1\leq i\leq n-1$, except possibly
for $i=k+1$ when $n=2k+2$.
\end{lemma}

\begin{proof}
Let $r\coloneqq\min\{k,n-1\}$. Projection of $E$ to $M$ has
contractible fibers $H_{\nu(p)}$, and projection of $E$ to $\Sph^n$
has $k$-connected fibers. So Lemma~\ref{lem:micro} yields
$\pi_i(M)\cong\pi_i(E)\cong\pi_i(\Sph^n)=0$ for $1\leq i\leq r$,
which proves the first assertion.
Thus $H_i(M;\Z)=0$ for $1\leq i\leq r$ by the Hurewicz theorem,
and then $H^i(M;\Z)=0$ for $1\leq i\leq r$ by the universal
coefficient theorem. So Poincar\'e duality yields $H_i(M;\Z)=0$ for
$n-r\leq i\leq n-1$. Now suppose that $n\leq 2k+2$. If $k\geq n-1$,
we are done. Otherwise $r=k$, and the only degrees not covered are
$k+1\leq i\leq n-k-1$. This range is empty if $n\leq 2k+1$, and
consists of $i=k+1$ if $n=2k+2$.
\end{proof}

\begin{lemma}\label{lem:complement}
Let $S\subset M$ be a nonempty proper open subset, and
$S'\coloneqq M\setminus S$. Suppose that $S$ is $k$-connected,
$S'$ has arbitrarily small $k$-connected open neighborhoods,
and $H_i(M;\Z)=0$ for $1\leq i\leq n-1$. If $n\leq 2k+2$, then
$S$ is contractible.
\end{lemma}

\begin{proof}
By Alexander duality \cite[Thm.~6.2.17]{spanier1966},
$H_j(M,S;\Z)$ is isomorphic to the direct limit of
$H^{n-j}(U;\Z)$ over open neighborhoods $U$ of $S'$. Since
$H_i(M;\Z)=0$ for $1\leq i\leq n-1$, the long exact sequence of
$(M,S)$ then shows that $\tilde H_i(S;\Z)$ is isomorphic to the direct
limit of $\tilde H^{n-i-1}(U;\Z)$, for $0\leq i\leq n-1$.
By hypothesis there is a cofinal family of such $U$ which are
$k$-connected, so that $\tilde H^{m}(U;\Z)=0$ for $m\leq k$.
Thus $\tilde H_i(S;\Z)=0$ for $n-k-1\leq i\leq n-1$, while
$\tilde H_i(S;\Z)=0$ for $i\leq k$ since $S$ is $k$-connected, and
$H_i(S;\Z)=0$ for $i\geq n$ since $S$ is a noncompact $n$-manifold.
As $n\leq 2k+2$, these ranges cover all degrees, so $S$ is acyclic.
Furthermore $S$ is simply connected if $k\geq1$, and also if $k=0$,
since then $n=2$, the fundamental group of a noncompact surface is
free, and $H_1(S;\Z)=0$ is its abelianization.
Hence $S$ is contractible by the Hurewicz and Whitehead theorems
\cite{hatcher2002}, since it has the homotopy type of a CW complex
\cite{milnor1959}.
\end{proof}

\begin{proposition}\label{prop:contractible-shades}
Suppose that every shade is $k$-connected, and either $n\leq 2k+1$,
or $n=2k+2$ and $H_{k+1}(M;\Z)=0$. Then every shade is contractible.
\end{proposition}

\begin{proof}
By Lemma~\ref{lem:homology}, $H_i(M;\Z)=0$ for $1\leq i\leq n-1$.
Now fix $u\in\Sph^n$. Then $S_u'=\{\sigma_u\leq0\}$ is nonempty
because $S_{-u}\neq\varnothing$, and its neighborhoods
$\{\sigma_u<c\}=S_{-u,c}$, $0<c<1$, shrink to $S_u'$ and are
$k$-connected by Lemma~\ref{lem:enlarged}. So
Lemma~\ref{lem:complement} shows that $S_u$ is contractible.
\end{proof}

\subsection{Proof of Theorem~\ref{thm:main}}
The first assertion is Corollary~\ref{cor:contractible}, and the
second follows from Corollary~\ref{cor:contractible} together with
Proposition~\ref{prop:contractible-shades}.\qed

\begin{note}\label{note:borderline}
When $k=0$, the hypothesis $H_{k+1}(M;\Z)=0$ in Theorem~\ref{thm:main}
cannot be dropped, by the counterexamples in \cite{ghomi2002-shadows,ghomi-tabachnikov2008}, which are constructed using ``skew loops" \cite{ghomi-solomon2002}.
For $k\geq1$, whether it can be dropped when $n=2k+2$ remains open.
In the latter case Lemma~\ref{lem:homology} shows that $M$ is $k$-connected, and
$H_{k+1}(M;\Z)$ is the only homology group of $M$ in degrees $1$
through $n-1$ which may not vanish. So
$H_{k+1}(M;\Z)\cong H^{k+1}(M;\Z)\cong\operatorname{Hom}(H_{k+1}(M;\Z),\Z)$
by Poincar\'e duality and the universal coefficient theorem; thus
$H_{k+1}(M;\Z)$ is free, and the hypothesis is equivalent to the
vanishing of the Betti number $b_{k+1}(M)$. When $k=0$, it is also
equivalent to $M\cong\Sph^2$, by the classification of surfaces.
When $k=1$, it is equivalent to $M\cong\Sph^4$, since $M$ is then a
simply connected homology $4$-sphere, hence a homotopy sphere, and so
is homeomorphic to $\Sph^4$ by Freedman's theorem \cite{freedman1982}.
The simplest candidate for a
counterexample is $M=\Sph^2\times\Sph^2$,
which is the union of two copies of $\Sph^2\times\B^2$ along
$\Sph^2\times\Sph^1$; we do not know whether it admits an embedding
in $\R^5$ with simply connected shades.
\end{note}

\section{Nonconvex Embeddings}\label{sec:joins}

In this section we prove Theorem~\ref{thm:counterexamples}. Our
examples will be parametrized by spherical joins
\cites{hatcher2002,bridson-haefliger1999}. For $p,q\geq1$, the
\emph{spherical join} of $\Sph^p$ and $\Sph^q$ is
\[
\Sph^p*\Sph^q\coloneqq
\bigl([0,\pi/2]\times\Sph^p\times\Sph^q\bigr)/\!\sim,
\]
where points of the product are denoted by $(t,x,y)$, and
$(0,x,y)\sim(0,x,y')$, $(\pi/2,x,y)\sim(\pi/2,x',y)$; so the $\Sph^q$
factor is collapsed at $t=0$, and the $\Sph^p$ factor at $t=\pi/2$.
The equivalence class of $(t,x,y)$ is denoted by $[t,x,y]$.
The standard map
$[t,x,y]\mapsto(\cos t\,x,\sin t\,y)$ identifies $\Sph^p*\Sph^q$
smoothly with the unit sphere of $\R^{p+1}\oplus\R^{q+1}$.
Theorem~\ref{thm:counterexamples} follows from the next result by
setting $p\coloneqq k+1$ and $q\coloneqq n-k-2$, which are at least
$k+1$ when $n\geq 2k+3$.

\begin{theorem}\label{thm:folded-join}
Let $p,q\geq1$. There is a smooth nonconvex
embedding $f\colon\Sph^p*\Sph^q\to\R^{p+q+2}$ whose shades are all
$(\min\{p,q\}-1)$-connected.
\end{theorem}

We construct $f$ and establish its nonconvexity in
Section~\ref{subsec:construction}. The topology of its shades is studied in Section~\ref{subsec:shades},
and the proof is completed in Section~\ref{subsec:proof}.

\subsection{The construction}\label{subsec:construction}
 Fix $p,q\geq1$, and put $n\coloneqq p+q+1$ and
$M\coloneqq\Sph^p*\Sph^q$. Near the endpoints of $[0,\pi/2]$ we use the
standard charts
\[
(x,z)\in\Sph^p\times\R^{q+1},\quad z=ty,
\qquad
(w,y)\in\R^{p+1}\times\Sph^q,\quad w=(\pi/2-t)x,
\]
where $|z|,|w|<\pi/2$.
In the above identification of $\Sph^p*\Sph^q$ with $\Sph^n$, the Gauss
map is $[t,x,y]\mapsto(\cos t\,x,\sin t\,y)$. The idea is to replace $t$
in this formula by a nonmonotone function $\theta(t)$, and then integrate
to find a hypersurface with the resulting Gauss map.
Let $\theta\colon\R\to\R$ be a smooth odd function such that
$\theta(\pi/2-t)=\pi/2-\theta(t)$, $0<\theta(t)<\pi/2$ for
$0<t<\pi/2$, and $\theta'(t_0)<0$ for some $t_0\in(0,\pi/2)$.
For instance we may set $\theta(t)\coloneqq t+\frac13\sin(4t)$. Also let
$a$, $b\colon[0,\pi/2]\to\R$ be the functions determined by
\[
a'(t)=-\sin\theta(t),\quad a(\pi/2)=0,
\qquad\quad
b'(t)=\cos\theta(t),\quad b(0)=0.
\]
Then $a'<0$ and $b'>0$ on $(0,\pi/2)$; in particular $a$ and $b$ are
positive on $(0,\pi/2)$. Define
\[
f([t,x,y])\coloneqq(a(t)x,\,b(t)y)\in\R^{p+1}\oplus\R^{q+1}.
\]

\begin{lemma}\label{lem:join-embedding}
$f$ is a smooth embedding with the Gauss map
\[
\nu([t,x,y])\coloneqq(\cos\theta(t)\,x,\sin\theta(t)\,y).
\]
\end{lemma}

\begin{proof}
For $0<t<\pi/2$, $f_t=(-\sin\theta(t)\,x,\cos\theta(t)\,y)$, while
$f_\xi=(a(t)\xi,0)$ and $f_\eta=(0,b(t)\eta)$ for $\xi\in T_x\Sph^p$ and
$\eta\in T_y\Sph^q$. These vectors are linearly independent because
$a(t),b(t)>0$ and $\norm{f_t}=1$. The vector $\nu([t,x,y])$ has unit
length and is orthogonal to all of them.
Endpoint smoothness follows from parity. Since $\theta$ is odd, there
are smooth one-variable functions
$\alpha,\beta,\kappa,\delta$ such that
\[
a(0)-a(t)=t^2\alpha(t^2),\qquad
b(t)=t\beta(t^2),\qquad
\cos\theta(t)=\kappa(t^2),\qquad
\sin\theta(t)=t\delta(t^2),
\]
with $\beta(0)=1$. In the chart
$z=ty$,
\[
f(x,z)=\bigl(\bigl(a(0)-\norm{z}^2\alpha(\norm{z}^2)\bigr)x,\;
\beta(\norm{z}^2)z\bigr),
\qquad
\nu(x,z)=\bigl(\kappa(\norm{z}^2)x,\;\delta(\norm{z}^2)z\bigr).
\]
Thus both maps extend smoothly across $z=0$, and
$df_{(x,0)}(\xi,v)=(a(0)\xi,v)$ has full rank.
At the other endpoint write $s=\pi/2-t$. The symmetry of $\theta$ gives
$a(\pi/2-s)=b(s)$ and $b(\pi/2-s)=a(s)$. The same parity calculation in
the chart $w=sx$ gives smooth extensions and
$df_{(0,y)}(v,\eta)=(v,b(\pi/2)\eta)$, again of full rank.
Finally, the norms of the two Euclidean blocks of $f([t,x,y])$ are $a(t)$
and $b(t)$. Their strict monotonicity determines $t$, after which every
noncollapsed spherical variable is determined. Thus $f$ is injective, and so is an embedding since $M$ is compact.
\end{proof}

\begin{lemma}\label{lem:join-nonconvex}
$f$ is not convex.
\end{lemma}

\begin{proof}
Let $\mathrm{II}(V,W)\coloneqq\ip{d\nu(V)}{df(W)}$ be the second fundamental form of $f$. On
$0<t<\pi/2$, $d\nu(\partial_t)=\theta'(t)f_t$, so
$\mathrm{II}(\partial_t,\partial_t)=\theta'(t)$. For $\xi\in T_x\Sph^p$
and $\eta\in T_y\Sph^q$,
\[
\mathrm{II}(\xi,\xi)=a(t)\cos\theta(t)\norm{\xi}^2,
\qquad
\mathrm{II}(\eta,\eta)=b(t)\sin\theta(t)\norm{\eta}^2.
\]
So the principal curvatures of $f$ are
$\theta'(t)$, $\cos\theta(t)/a(t)$, $\sin\theta(t)/b(t)$.
At $t=t_0$, the last two are positive while $\theta'(t_0)<0$.
So $\mathrm{II}$ is indefinite. 
\end{proof}

\subsection{Topology of the shades}\label{subsec:shades}
Fix $u=(u_1,u_2)\in\R^{p+1}\oplus\R^{q+1}$,
so that $|u_1|^2+|u_2|^2=1$,
and choose $e_1\in\Sph^p$ and $e_2\in\Sph^q$
so that $u_i=|u_i|\,e_i$;
if $u_i=0$, choose $e_i$ arbitrarily.
Put $x_1\coloneqq\ip{x}{e_1}$ and $y_2\coloneqq\ip{y}{e_2}$.
Then the shading function is
\begin{equation}\label{eq:shade-H}
\sigma_u([t,x,y])=|u_1|\cos\theta(t)\,x_1+|u_2|\sin\theta(t)\,y_2.
\end{equation}

\begin{proposition}\label{prop:shade-simply}
The shade $S_u$ is $(\min\{p,q\}-1)$-connected.
\end{proposition}

\begin{proof}
Let $0\leq j<\min\{p,q\}$, and $g\colon\Sph^j\to S_u$ be a continuous
map. The spheres $\{t=0\}$ and $\{t=\pi/2\}$ have codimensions
$q+1$ and $p+1$ in $M$, and in their complement
$\{0<t<\pi/2\}\cong(0,\pi/2)\times\Sph^p\times\Sph^q$
the sets $\{x=-e_1\}$ and $\{y=-e_2\}$ have codimensions
$p$ and $q$. All these exceed $j$.
So by transversality, we may assume, after a small homotopy
in $S_u$, that $g$ avoids these sets, i.e., $g=[t,x,y]$ with
$0<t<\pi/2$, $x\neq-e_1$, and $y\neq-e_2$.
Then the shortest great circle arcs from $x$ to $e_1$, and from $y$
to $e_2$, depend continuously on the point of $\Sph^j$. Moving $x$
and $y$ along these arcs, keeping $t$ fixed, does not decrease $x_1$
or $y_2$, and so by \eqref{eq:shade-H} does not decrease $\sigma_u$.
Thus $g$ is homotopic in $S_u$ to a map into the central ridge
$\{[t,e_1,e_2]:0<t<\pi/2\}$, where
$\sigma_u=|u_1|\cos\theta(t)+|u_2|\sin\theta(t)>0$. Finally, moving
$t$ linearly to $\pi/4$ contracts this map within the ridge. Hence
$\pi_j(S_u)=0$.
\end{proof}

\subsection{Proof of Theorem~\ref{thm:folded-join}}\label{subsec:proof}

Let $f$ be the map constructed in Section~\ref{subsec:construction}.
Then $f$ is a smooth embedding by Lemma~\ref{lem:join-embedding},
which is not convex by Lemma~\ref{lem:join-nonconvex}, and its shades
are $(\min\{p,q\}-1)$-connected by
Proposition~\ref{prop:shade-simply}.\qed

\section*{Acknowledgement}
The author thanks Shibing Chen and Yuanyuan Li for their
interest in this problem, and
Igor Belegradek for helpful comments.
AI tools (Claude by Anthropic and ChatGPT by OpenAI) were used
in preparing and editing the manuscript. The author has
verified all AI-assisted content and takes full responsibility
for the final text.

\bibliography{references}

\begin{bibdiv}
\begin{biblist}

\bib{bridson-haefliger1999}{book}{
      author={Bridson, Martin~R.},
      author={Haefliger, Andr\'{e}},
       title={Metric spaces of non-positive curvature},
      series={Grundlehren der mathematischen Wissenschaften [Fundamental
  Principles of Mathematical Sciences]},
   publisher={Springer-Verlag, Berlin},
        date={1999},
      volume={319},
        ISBN={3-540-64324-9},
         url={https://doi.org/10.1007/978-3-662-12494-9},
      review={\MR{1744486}},
}

\bib{chern-lashof1957}{article}{
      author={Chern, Shiing-Shen},
      author={Lashof, Richard~K.},
       title={On the total curvature of immersed manifolds},
        date={1957},
        ISSN={0002-9327},
     journal={Amer. J. Math.},
      volume={79},
       pages={306\ndash 318},
         url={http://dx.doi.org.prx.library.gatech.edu/10.2307/2372684},
      review={\MR{0084811}},
}

\bib{chern-lashof1958}{article}{
      author={Chern, Shiing-Shen},
      author={Lashof, Richard~K.},
       title={On the total curvature of immersed manifolds. {II}},
        date={1958},
        ISSN={0026-2285},
     journal={Michigan Math. J.},
      volume={5},
       pages={5\ndash 12},
  url={http://projecteuclid.org.prx.library.gatech.edu/euclid.mmj/1028998005},
      review={\MR{0097834}},
}

\bib{choe1990}{article}{
      author={Choe, Jaigyoung},
       title={Index, vision number and stability of complete minimal surfaces},
        date={1990},
        ISSN={0003-9527},
     journal={Arch. Rational Mech. Anal.},
      volume={109},
      number={3},
       pages={195\ndash 212},
      review={\MR{1025170}},
}

\bib{freedman1982}{article}{
      author={Freedman, Michael~Hartley},
       title={The topology of four-dimensional manifolds},
        date={1982},
        ISSN={0022-040X,1945-743X},
     journal={J. Differential Geometry},
      volume={17},
      number={3},
       pages={357\ndash 453},
         url={http://projecteuclid.org/euclid.jdg/1214437136},
      review={\MR{679066}},
}

\bib{ghomi2002-shadows}{article}{
      author={Ghomi, Mohammad},
       title={Shadows and convexity of surfaces},
        date={2002},
        ISSN={0003-486X,1939-8980},
     journal={Ann. of Math. (2)},
      volume={155},
      number={1},
       pages={281\ndash 293},
         url={https://doi.org/10.2307/3062157},
      review={\MR{1888801}},
}

\bib{ghomi2002-solution}{incollection}{
      author={Ghomi, Mohammad},
       title={Solution to the shadow problem in 3-space},
        date={2002},
   booktitle={Minimal surfaces, geometric analysis and symplectic geometry
  ({B}altimore, {MD}, 1999)},
      series={Adv. Stud. Pure Math.},
      volume={34},
   publisher={Math. Soc. Japan, Tokyo},
       pages={129\ndash 142},
         url={https://doi.org/10.2969/aspm/03410129},
      review={\MR{1925735}},
}

\bib{ghomi2007}{article}{
      author={Ghomi, Mohammad},
       title={Topology of surfaces with connected shades},
        date={2007},
        ISSN={1093-6106},
     journal={Asian J. Math.},
      volume={11},
      number={4},
       pages={621\ndash 634},
      review={\MR{2402941}},
}

\bib{ghomi-solomon2002}{article}{
      author={Ghomi, Mohammad},
      author={Solomon, Bruce},
       title={Skew loops and quadric surfaces},
        date={2002},
        ISSN={0010-2571},
     journal={Comment. Math. Helv.},
      volume={77},
      number={4},
       pages={767\ndash 782},
      review={\MR{1949113}},
}

\bib{ghomi-tabachnikov2008}{article}{
      author={Ghomi, Mohammad},
      author={Tabachnikov, Serge},
       title={Totally skew embeddings of manifolds},
        date={2008},
        ISSN={0025-5874},
     journal={Math. Z.},
      volume={258},
      number={3},
       pages={499\ndash 512},
         url={http://dx.doi.org/10.1007/s00209-007-0182-8},
      review={\MR{2369041}},
}

\bib{hatcher2002}{book}{
      author={Hatcher, Allen},
       title={Algebraic topology},
   publisher={Cambridge University Press},
     address={Cambridge},
        date={2002},
      review={\MR{1867354}},
}

\bib{kuiper1970}{article}{
      author={Kuiper, Nicolaas~H.},
       title={Minimal total absolute curvature for immersions},
        date={1970},
        ISSN={0020-9910},
     journal={Invent. Math.},
      volume={10},
       pages={209\ndash 238},
         url={http://dx.doi.org.prx.library.gatech.edu/10.1007/BF01403250},
      review={\MR{0267597}},
}

\bib{kuiper:tight}{incollection}{
      author={Kuiper, Nicolaas~H.},
       title={Geometry in curvature theory},
        date={1997},
   booktitle={Tight and taut submanifolds ({B}erkeley, {CA}, 1994)},
   publisher={Cambridge Univ. Press},
     address={Cambridge},
       pages={1\ndash 50},
      review={\MR{1486869}},
}

\bib{milnor1959}{article}{
      author={Milnor, John},
       title={On spaces having the homotopy type of a {${\rm CW}$}-complex},
        date={1959},
        ISSN={0002-9947,1088-6850},
     journal={Trans. Amer. Math. Soc.},
      volume={90},
       pages={272\ndash 280},
         url={https://doi.org/10.2307/1993204},
      review={\MR{100267}},
}

\bib{raptis2017}{article}{
      author={Raptis, George},
       title={On {S}erre microfibrations and a lemma of {M}. {W}eiss},
        date={2017},
        ISSN={0017-0895,1469-509X},
     journal={Glasg. Math. J.},
      volume={59},
      number={3},
       pages={649\ndash 657},
         url={https://doi.org/10.1017/S0017089516000458},
      review={\MR{3682004}},
}

\bib{spanier1966}{book}{
      author={Spanier, Edwin~H.},
       title={Algebraic topology},
   publisher={McGraw-Hill Book Co.},
     address={New York-Toronto-London},
        date={1966},
      review={\MR{0210112}},
}

\bib{weiss2005}{article}{
      author={Weiss, Michael},
       title={What does the classifying space of a category classify?},
        date={2005},
        ISSN={1532-0081,1532-0073},
     journal={Homology Homotopy Appl.},
      volume={7},
      number={1},
       pages={185\ndash 195},
         url={http://projecteuclid.org/euclid.hha/1139839512},
      review={\MR{2175298}},
}

\end{biblist}
\end{bibdiv}

\end{document}